\documentclass[10pt]{article}
\usepackage{amsmath, amssymb, amscd, amsthm}

\newcommand{\xra}[1]{\ensuremath{\xrightarrow{#1}}}
\newcommand{\xla}[1]{\ensuremath{\xleftarrow{#1}}}
\newcommand{\B}[1]{\ensuremath{\mathbb{#1}}}

\newcommand{\sgn}[1]{\ensuremath{\text{sgn}{#1}}}
\newcommand{\lcomod}[1]{\ensuremath{{#1}\text{--{\bf Comod}}}}

\newcommand{\eotimes}[1]{\ensuremath{\underset{#1}{\otimes}}}

\newtheorem{thm}{Theorem}[section]
\newtheorem{cor}[thm]{Corollary}
\newtheorem{lem}[thm]{Lemma}

\theoremstyle{definition}
\newtheorem{defn}[thm]{Definition}
\newtheorem{rem}[thm]{Remark}

\numberwithin{equation}{section}

\title{Loday--Quillen--Tsygan Theorem for Coalgebras}
\author{Atabey Kaygun}
\begin{document}
\maketitle

\section{Introduction}

The original Loday--Quillen--Tsygan Theorem (LQT) is proven by Loday and
Quillen \cite{LodayQuillen:LQT} and independently by Tsygan
\cite{Tsygan:LQT}.  It states that the ordinary Lie homology (here
referred as Chevalley--Eilenberg--Lie homology) of the Lie algebra of
the infinite matrices $gl(A)$ over an unital associative algebra $A$ is
generated by the cyclic homology of $A$ as an exterior algebra.
Although Lie algebras have been studied extensively, their
non-commutative counterparts, Leibniz algebras and their homology, are
defined only recently by Loday \cite{Loday:LeibnizAlgebras}.  In this
setting, LQT is extended by Cuvier \cite{Cuvier:LQT} proving Leibniz
homology (here referred as Chevalley--Eilenberg--Leibniz homology) of
$gl(A)$ is generated by the Hochschild homology of $A$ as a tensor
algebra.  There is also a slightly different proof of Cuvier's result by
Oudom \cite{Oudom:Leibniz} by using a specific filtration on the
Chevalley--Eilenberg--Leibniz complex of $gl(A)$ and a spectral
sequence.  All of these proofs and the proof we present here for the
coalgebras rely heavily on Weyl's invariant theory \cite[Chapter
9]{Loday:CyclicHomology}.  We would like to mention that in
\cite{AboughaziOgle:CyclicHomology} Aboughazi and Ogle gave an
alternative proof of LQT which did not use Weyl's invariant theory.
There is also a general LQT type result for algebras over operads by
Fresse \cite{Fresse:Cogroups}.

Coassociative coalgebras and their homologies received some attention
recently in the context of Hopf and bialgebra cyclic (co)homology
\cite{ConnesMoscovici:HopfCyclicCohomology,
Crainic:CyclicCohomologyOfHopfAlgebras, Khalkhali:HopfCyclicHomology,
Kaygun:BialgebraCyclic}.  On the other hand, Lie coalgebras have not
received much attention in their own right
\cite{Michaelis:LieCoalgebras, Oziewicz:LieCoalgebras,
Cuartero:LieCoalgebras} even though they appear as auxillary structures
in important results such as Hinich's explanation
\cite{Hinich:TamarkinKontsevic} of Tamarkin's proof of Kontsevic
Formality Theorem \cite{Tamarkin:FormalityI}.  The only reference to
Leibniz coalgebras we found was in Livernet
\cite{Livernet:RationalHomotopy}.

In this paper we prove that LQT generalizes to the case of coalgebras.
Specifically, we show that the Chevalley--Eilenberg--Lie homology of the
Lie coalgebra of matrices $gl^c(C)$ over a coassociative coalgebra $C$
is generated by the cyclic homology of the underlying coalgebra $C$ as
an exterior algebra.

Here is the plan of this paper.  In Section~\ref{Coalgebras} we give a
self-contained account of coassociative, Leibniz and Lie coalgebras and
their comodules.  In Section~\ref{Bar} and Section~\ref{Hochschild} we
develop the bar and Hochschild homology theories for coassociative
coalgebras.  In Section~\ref{CEL} we define several homology theories
for Leibniz coalgebras including Chevalley--Eilenberg--Leibniz homology
and symmetric Chevalley--Eilenberg--Leibniz homology which we call as
Chevalley--Eilenberg--Lie homology.  In Section~\ref{Tale} we connect
all of these homology theories together for the Lie coalgebra $Lie(C)$
of a coassociative coalgebra $C$.  Section~\ref{Misc} contains several
results about matrix coalgebras.  Finally, in Section~\ref{LQT} we prove
Loday--Quillen--Tsygan Theorem for coassociative counital coalgebras.

\noindent {\bf Acknowledgements. } We would like to thank Henri
Moscovici for his helpful remarks.

Throughout this paper, we will assume that $k$ is a field of
characteristic $0$.  The tensor products are taken over $k$ and are
denoted by $\otimes$.

\section{Coalgebras and their comodules}\label{Coalgebras}

\begin{defn}
A $k$--module $C$ is called a coassociative coalgebra if there is
morphism of $k$--modules $C\xra{\Delta}C\otimes C$ such that
\begin{align}
(id_C\otimes\Delta)\Delta 
  = & \ (\Delta\otimes id_C)\Delta
\end{align}
$C$ is called counital if there is a morphism of $k$--modules 
$C\xra{\eta}k$ such that 
\begin{align}
(\eta\otimes id_C)\Delta = & \ id_C =
(id_C\otimes\eta)\Delta
\end{align}
We will use Sweedler's notation and denote $\Delta(c)$ by $\sum_c
c_{(1)}\otimes c_{(2)}$ for any $c\in C$.
\end{defn}

\begin{defn}
Assume $(C,\Delta,\eta)$ is a counital coassociative coalgebra.  A
$k$--module $X$ is called a left $C$--comodule iff there is a morphism
of $k$--modules $X\xra{\rho_X}C\otimes X$ such that
\begin{align}
(id_C\otimes\rho_X)\rho_X    & = (\Delta\otimes id_X)\rho_X\\
(\eta\otimes id_C)\rho_X & = id_X
\end{align}
Right $C$--comodules and $C$--bicomodules are defined similarly.  Again,
we will use Sweedler's notation and denote $\rho_X(x)$ by $\sum_x
x_{(-1)}\otimes x_{(0)}$ for any $x\in X$ whenever $X$ is a left
$C$--comodule.
\end{defn}

\begin{defn}
A $k$--module $L$ is called a Leibniz coalgebra iff there is a morphism
of $k$--modules $L\xra{\delta}L\otimes L$ such that 
\begin{equation}
(id_3-(id_1\otimes\tau_2))(\delta\otimes id_1)\delta
=(id_1\otimes\delta)\delta
\end{equation}
where $L^{\otimes 2}\xra{\tau_2}L^{\otimes 2}$ is the permutation of the
factors. We will use Sweedler's notation and denote $\delta(x)$ by
$x_{[1]}\otimes x_{[2]}$ for any $x\in L$.  Then, for any $x\in L$ we
have
\begin{equation}\label{coJacobi}
\sum x_{[1][1]}\otimes x_{[1][2]}\otimes x_{[2]}
    -x_{[1][1]}\otimes x_{[2]}\otimes x_{[1][2]}
=\sum x_{[1]}\otimes x_{[2][1]}\otimes x_{[2][2]}
\end{equation}
The morphism $\delta$ is called a cobracket.
\end{defn}

\begin{defn}
A Leibniz coalgebra $(L,\delta)$ is called a Lie coalgebra iff $L$ is a
Leibniz coalgebra and the cobracket $\delta$ is anti-symmetric, i.e.
\begin{equation}
(id_2+\tau_2)\delta = 0
\end{equation}
\end{defn}

\begin{lem}
Let $(C,\Delta)$ be a coassociative coalgebra.  Let
$\delta=(id_2-\tau_2)\Delta$.  Then $(C,\delta)$ is a Lie coalgebra.
This Lie coalgebra associated with a coassociative coalgebra
$(C,\Delta)$ is denoted by $(Lie(C),\delta)$.
\end{lem}

\begin{proof}
It is easy to see that
\begin{equation*}
(id_2+\tau_2)\delta(x)=(id_2+\tau_2)(1-\tau_2)\Delta(x)
= (1-\tau_2^2)\Delta(x)=0
\end{equation*}
For the coJacobi condition coassociativity will play an important role.
Since $\Delta$ is coassociative, we can write
\begin{equation*}
\delta(x)= x_{(1)}\otimes x_{(2)}-x_{(2)}\otimes x_{(1)}
\end{equation*}
Then 
\begin{align*}
(id_3-(id_1\otimes\tau_2)) & (\delta\otimes id_1)\delta(x)\\
= & (id_3-(id_1\otimes\tau_2))(\delta\otimes id_1)\left(  
	 x_{(1)}\otimes x_{(2)}
       - x_{(2)}\otimes x_{(1)}\right)\\
= & (id_3-(id_1\otimes\tau_2))\left(
           x_{(1)(1)}\otimes x_{(1)(2)}\otimes x_{(2)}
         - x_{(1)(2)}\otimes x_{(1)(1)}\otimes x_{(2)}\right)\\
& + (id_3-(id_1\otimes\tau_2))\left(
         - x_{(2)(1)}\otimes x_{(2)(2)}\otimes x_{(1)}
         + x_{(2)(2)}\otimes x_{(2)(1)}\otimes x_{(1)}\right)\\
= & (id_3-(id_1\otimes\tau_2))\left(
           x_{(1)}\otimes x_{(2)}\otimes x_{(3)}
         - x_{(2)}\otimes x_{(1)}\otimes x_{(3)} \right)\\
& + (id_3-(id_1\otimes\tau_2))\left(
         - x_{(2)}\otimes x_{(3)}\otimes x_{(1)}
         + x_{(3)}\otimes x_{(2)}\otimes x_{(1)}\right)\\
= &        x_{(1)}\otimes x_{(2)}\otimes x_{(3)}
         - x_{(1)}\otimes x_{(3)}\otimes x_{(2)}
         - x_{(2)}\otimes x_{(1)}\otimes x_{(3)}
         + x_{(2)}\otimes x_{(3)}\otimes x_{(1)}\\
&        - x_{(2)}\otimes x_{(3)}\otimes x_{(1)}  
         + x_{(2)}\otimes x_{(1)}\otimes x_{(3)}
         + x_{(3)}\otimes x_{(2)}\otimes x_{(1)}
         - x_{(3)}\otimes x_{(1)}\otimes x_{(2)}
\end{align*}
and after the cancellations we get
\begin{align*}
(id_3-(id_1\otimes\tau_2))(\delta\otimes id_1)\delta(x)
= &        x_{(1)}\otimes x_{(2)}\otimes x_{(3)}
         - x_{(1)}\otimes x_{(3)}\otimes x_{(2)}\\
  &      + x_{(3)}\otimes x_{(2)}\otimes x_{(1)}
         - x_{(3)}\otimes x_{(1)}\otimes x_{(2)}\\
= &        x_{(1)}\otimes x_{(2)(1)}\otimes x_{(2)(2)}
         - x_{(1)}\otimes x_{(2)(2)}\otimes x_{(2)(1)}\\
  &      + x_{(2)}\otimes x_{(1)(1)}\otimes x_{(1)(2)}
         - x_{(2)}\otimes x_{(1)(2)}\otimes x_{(1)(2)}\\
= & x_{(1)}\otimes\delta(x_{(2)}) + x_{(2)}\otimes\delta(x_{(1)})\\
= & (id_1\otimes\delta)\delta(x)
\end{align*}
as we wanted to show.
\end{proof}

\begin{defn}
Let $L$ be a Leibniz coalgebra.  A $k$--module $X$ is called a right
$L$--comodule iff there is a morphism of $k$--modules
$X\xra{\rho_X}X\otimes L$ such that
\begin{align}
x_{[0][0]}\otimes x_{[0][1]}\otimes x_{[1]}
 - x_{[0][0]}\otimes x_{[1]}\otimes x_{[0][1]}
 = x_{[0]}\otimes x_{[1][1]}\otimes x_{[1][2]}
\end{align}
for any $x\in X$ where we write $\rho_X(x)=x_{[0]}\otimes x_{[1]}$
\end{defn}

\begin{lem}\label{DiagonalCoaction}
Let $X$ and $Y$ be two right $L$--comodule.  Then $X\otimes Y$ has a
right $L$--comodule structure.
\end{lem}

\begin{proof}
Let $(X\otimes Y)\xra{\rho_{X\otimes Y}}(X\otimes Y)\otimes L$ be
defined as 
\begin{align}
\rho_{X\otimes Y}(x\otimes y)
 = (x_{[0]}\otimes y)\otimes x_{[1]} + (x\otimes y_{[0]})\otimes y_{[1]}
\end{align}
for any $(x\otimes y)$ from $X\otimes Y$.  Then
\begin{align*}
(x\otimes y)_{[0][0]}&\otimes(x\otimes y)_{[0][1]}\otimes(x\otimes
  y)_{[1]} 
-  (x\otimes y)_{[0][0]}\otimes(x\otimes y)_{[1]}\otimes(x\otimes
  y)_{[0][1]} \\
 = & (x_{[0][0]}\otimes y)\otimes x_{[0][1]}\otimes x_{[1]}
     + (x_{[0]}\otimes y_{[0]})\otimes x_{[1]}\otimes y_{[1]}
     + (x_{[0]}\otimes y_{[0]})\otimes y_{[1]}\otimes x_{[1]}\\ 
   & + (x\otimes y_{[0][0]})\otimes y_{[0][1]}\otimes y_{[1]}
     - (x_{[0][0]}\otimes y)\otimes x_{[1]}\otimes x_{[0][1]}
     - (x_{[0]}\otimes y_{[0]})\otimes y_{[1]}\otimes x_{[1]}\\
   & - (x_{[0]}\otimes y_{[0]})\otimes x_{[1]}\otimes y_{[1]} 
     - (x\otimes y_{[0][0]})\otimes y_{[1]}\otimes y_{[0][1]}\\
 = & (x_{[0][0]}\otimes y)\otimes x_{[0][1]}\otimes x_{[1]}
     - (x_{[0][0]}\otimes y)\otimes x_{[1]}\otimes x_{[0][1]}\\
   & + (x\otimes y_{[0][0]})\otimes y_{[0][1]}\otimes y_{[1]}
     - (x\otimes y_{[0][0]})\otimes y_{[1]}\otimes y_{[0][1]}\\
 = & (x_{[0]}\otimes y)\otimes x_{[1][1]}\otimes x_{[1][2]}
    + (x\otimes y_{[0]})\otimes y_{[1][1]}\otimes y_{[1][2]}\\
 = & (x\otimes y)_{[0]}\otimes (x\otimes y)_{[1][1]}\otimes 
     (x\otimes y)_{[1][2]}
\end{align*}
as we wanted to show.
\end{proof}

\begin{rem}
For a Leibniz coalgebra $L$, kernels and cokernels in the category of
$L$--comodules are defined as kernels and cokernels in the category of
$k$--modules.  Then the category of $L$--comodules is an abelian
subcategory of the category of $k$--modules.
\end{rem}

\begin{thm}
Let $(C,\Delta)$ be a coassociative coalgebra.  Then there is a functor
from the category of left $C$--comodules into the category of left
$Lie(C)$--comodules of the form $\lcomod{C}\xra{Lie}\lcomod{Lie(C)}$.
\end{thm}

\begin{proof}
The functor $Lie$ is identity on both object and morphisms.  In order to
show that $Lie(X)$ is a $Lie(C)$--comodule whenever
$X\xra{\rho_X}C\otimes X$ is a $C$--comodule.  We need to define a
$Lie(C)$--comodule structure $X\xra{\rho_X}Lie(C)\otimes X$.  So, take
the same $\rho_X$.  Then
\begin{align*}
x_{[0][0]}\otimes x_{[0][1]}\otimes x_{[1]}
 - x_{[0][0]}\otimes x_{[1]}\otimes x_{[0][1]}
= & x_{(0)}\otimes x_{(1)}\otimes x_{(2)}
 - x_{(0)}\otimes x_{(2)}\otimes x_{(1)}\\
= & x_{(0)}\otimes x_{(1)[1]}\otimes x_{(1)[2]}\\
= & x_{[0]}\otimes x_{[1][1]}\otimes x_{[1][2]}
\end{align*}
which proves $\rho_X$ defines $Lie(C)$--comodule structure on $X$.
\end{proof}

\begin{defn}
Let $L$ be a Leibniz coalgebra and $M\xra{\rho_M}M\otimes L$ be a right
$L$--comodule.  Then $M^L$ the subcomodule of $C$ invariants of $M$ is
defined as
\begin{align*}
M^L = ker(M\xra{\rho_M}M\otimes L) 
\end{align*}
\end{defn}

\section{The bar complex}\label{Bar}

From this section on, assume $C$ is a counital coassociative coalgebra.
Let $Lie(C)$ be the Lie coalgebra associated with the coassociative
coalgebra $C$.

\begin{defn}
Let $CB_*(C):=\left\{C^{\otimes n}\right\}_{n\geq 0}$ and define
$d^{CB}_0 = 0$ and also let
\begin{align}
d^{CB}_n &= \sum_{j=0}^{n-1}(-1)^j(id_j\otimes\Delta\otimes id_{n-1-j})
\end{align}
for any $n\geq 1$.
\end{defn}

\begin{thm}
$CB_*(C)$ is a differential graded algebra.
\end{thm}

\begin{proof}
For $0\leq j\leq n-1$ define morphisms $\partial_j$ as
\begin{align}
\partial_j(c^1\otimes\cdots\otimes c^n)
 = (id_j\otimes\Delta\otimes id_{n-j})
   (c^1\otimes\cdots\otimes c^n)
 = (\cdots\otimes c^{j+1}_{(1)}\otimes c^{j+1}_{(2)}\otimes\cdots)
\end{align}
One can check that, for $i< j$
\begin{align*}
\partial_j\partial_i(c^0\otimes\cdots\otimes c^n)
 = & \begin{cases}
     (\cdots\otimes c^{i+1}_{(1)}\otimes c^{i+1}_{(2)}\otimes
      \cdots\otimes c^j_{(1)}\otimes c^j_{(2)}\otimes
      \cdots)         & \text{ if } i+1<j\\
     (\cdots\otimes c^{i+1}_{(1)}\otimes c^{i+1}_{(2)}\otimes 
      c^{i+1}_{(3)}\otimes\cdots)  & \text{ if } i+1=j
     \end{cases}\\
 = &\ \partial_i\partial_{j-1}(c^0\otimes\cdots\otimes c^n)
\end{align*}
Therefore $d^{CB}_n = \sum_{j=0}^{n-1}(-1)^j\partial_j$ is a
differential.   Note that
\begin{align}
d^{CB}_{n+m}
 = & \sum_{j=1}^n(-1)^{j-1}(id_{j-1}\otimes\Delta\otimes
         id_{n-j})\otimes id_m
     + \sum_{j=1}^m(-1)^{n+j-1}id_n\otimes(id_{j-1}\otimes\Delta\otimes 
         id_{m-j})\\
 = & (d^{CB}_n\otimes id_m) + (-1)^n(id_n\otimes d^{CB}_m)
\end{align}
And if we take the tensor multiplication as the algebra structure on
$CB_*(C)$, then with this equality
\begin{align*}
d^{CB}_{n+m}(\Psi\otimes\Phi)
 = (d^{CB}_n(\Psi)\otimes\Phi) + (-1)^n(\Psi\otimes d^{CB}_m(\Phi))
\end{align*}
for any $\Psi$ from $CB_n(C)$ and $\Phi$ from $CB_m(C)$ where $n\geq 0$
and $m\geq 0$ are arbitrary,  This proves that $CB_*(C)$ is a
differential graded $k$--algebra.
\end{proof}

\begin{thm}\label{BarComplex}
$CB_*(C)$ is a differential graded $Lie(C)$--comodule.
\end{thm}

\begin{proof}
The comodule structure is given by the diagonal coaction
\begin{align}
\rho_n(x^1\otimes\cdots\otimes x^n)
 = & \sum_{j=1}^n(\cdots\otimes x^j_{[1]}\otimes\cdots)
                 \otimes x^j_{[2]}\nonumber\\
 = & \sum_{j=1}^n(\cdots\otimes x^j_{(1)}\otimes\cdots)
                 \otimes x^j_{(2)}
    - \sum_{j=1}^n(\cdots\otimes x^j_{(2)}\otimes\cdots)
                 \otimes x^j_{(1)}
\end{align}
Lemma~\ref{DiagonalCoaction} proves that $CB_*(C)$ is a graded
$Lie(C)$--comodule.  In order to show that $CB_*(C)$ is a differential
graded $Lie(C)$--comodule, we must show that the coaction and the
differentials commute.  In other words we must show $\rho_{n+1}d^{CB}_n
= (d^{CB}_n\otimes id_1)\rho_n$ for any $n\geq 0$.  So, consider
\begin{align*}
\rho_{n+1}d^{CB}_n(x^1\otimes\cdots\otimes x^n)
 = & \sum_{i=0}^{n-1}\sum_{j=1}^i(-1)^i(\cdots\otimes x^j_{[1]}\otimes
                        \cdots\otimes x^{i+1}_{(1)}\otimes
                        x^{i+1}_{(2)}\otimes\cdots)\otimes x^j_{[2]}\\
   & + \sum_{i=0}^{n-1}(-1)^i(\cdots\otimes x^{i+1}_{(1)[1]}\otimes
                        x^{i+1}_{(2)}\otimes\cdots)\otimes x^{i+1}_{(1)[2]}\\
   & + \sum_{i=0}^{n-1}(-1)^i(\cdots\otimes x^{i+1}_{(1)}\otimes
                        x^{i+1}_{(2)[1]}\otimes\cdots)\otimes
                        x^{i+1}_{(2)[2]}\\
   & + \sum_{i=0}^{n-1}\sum_{j=i+3}^n(-1)^i(\cdots\otimes x^{i+1}_{(1)}
                        \otimes x^{i+1}_{(2)}\otimes\cdots\otimes x^{j-1}_{[1]}
                        \otimes\cdots)\otimes x^{j-1}_{[2]}
\end{align*}
However,
\begin{align*}
(x_{(1)[1]}&\otimes x_{(2)}\otimes x_{(1)[2]})
 + (x_{(1)}\otimes x_{(2)[1]}\otimes x_{(2)[2]})\\
 = &   (x_{(1)}\otimes x_{(3)}\otimes x_{(2)})
     - (x_{(2)}\otimes x_{(3)}\otimes x_{(1)})
     + (x_{(1)}\otimes x_{(2)}\otimes x_{(3)})
     - (x_{(1)}\otimes x_{(3)}\otimes x_{(2)})\\
 = &   (x_{(1)}\otimes x_{(2)}\otimes x_{(3)})
     - (x_{(2)}\otimes x_{(3)}\otimes x_{(1)})\\
 = &   (x_{(1)(1)}\otimes x_{(1)(2)}\otimes x_{(2)})
     - (x_{(2)(1)}\otimes x_{(2)(2)}\otimes x_{(1)})\\
 = &   (x_{[1](1)}\otimes x_{[1](2)}\otimes x_{[2]})
\end{align*}
for any $x\in C$.  Therefore
\begin{align*}
\rho_{n+1}d^{CB}_n (x^1\otimes\cdots\otimes x^n)
 = &   \sum_{i=0}^{n-1}\sum_{j=i+2}^{n-1}(-1)^i(\cdots\otimes x^{i+1}_{(1)}
                        \otimes x^{i+1}_{(2)}\otimes\cdots\otimes x^j_{[1]}
                        \otimes \cdots)\otimes x^j_{[2]}\\
   & + \sum_{i=0}^{n-1}(-1)^i(\cdots\otimes x^{i+1}_{[1](1)}\otimes
                        x^{i+1}_{[1](2)}\otimes\cdots)\otimes x^{i+1}_{[2]}\\
   & + \sum_{i=0}^{n-1}\sum_{j=1}^i(-1)^i(\cdots\otimes x^j_{[1]}\otimes
                        \cdots\otimes x^{i+1}_{(1)}\otimes
                        x^{i+1}_{(2)}\otimes\cdots)\otimes x^j_{[2]}\\
 = &  \sum_{j=1}^{n-1}\sum_{i=0}^{j-2}(-1)^i(\cdots\otimes x^{i+1}_{(1)}
                        \otimes x^{i+1}_{(2)}\otimes\cdots\otimes x^j_{[1]}
                        \otimes \cdots)\otimes x^j_{[2]}\\
   & + \sum_{i=0}^{n-1}(-1)^i(\cdots\otimes x^{i+1}_{[1](1)}\otimes
                        x^{i+1}_{[1](2)}\otimes\cdots)\otimes x^{i+1}_{[2]}\\
   & + \sum_{j=1}^{n-1}\sum_{i=j}^{n-1}(-1)^i(\cdots\otimes x^j_{[1]}\otimes
                        \cdots\otimes x^{i+1}_{(1)}\otimes
                        x^{i+1}_{(2)}\otimes\cdots)\otimes x^j_{[2]}\\
 = & (d^{CB}_n\otimes id_1)\rho_n(x^1\otimes\cdots\otimes x^n)
\end{align*}
as we wanted to show.
\end{proof}

\begin{thm}\label{BarHomotopy}
The $Lie(C)$--comodule structure on $CB_*(C)$ is null-homotopic.
\end{thm}

\begin{proof}
We need to provide a null-homotopy, i.e. a morphism of graded modules of
the form $CB_*(C)\xra{h_*}CB_*(C)[-1]\otimes Lie(C)$ which satisfies
$(d^{CB}_{n-1}\otimes id_C)h_n + h_{n+1}d^{CB}_n=\rho_n$.  Define
\begin{align}
h_n(x^1\otimes\cdots\otimes x^n)
 = & \sum_{j=1}^{n}(-1)^{j}(\cdots\otimes\widehat{x^{j}}\otimes\cdots)
                 \otimes x^{j}
\end{align}
for any $(x^1\otimes\cdots\otimes x^n)$ from $CB_n(C)$.  Then
\begin{align*}
(-1)^ih_{n+1}\partial_i(x^1\otimes\cdots\otimes x^n)
 = & \sum_{j=0}^{i-1}(-1)^{i+j+1}(\cdots\otimes\widehat{x^{j+1}}\otimes
                \cdots\otimes x^{i+1}_{(1)}\otimes x^{i+1}_{(2)}\otimes
                \cdots)\otimes x^{j+1}\\
   & + (\cdots\otimes x^{i+1}_{[1]}\otimes\cdots)\otimes x^{i+1}_{[2]}\\
   & + \sum_{j=i+2}^n(-1)^{i+j+1}(\cdots\otimes x^{i+1}_{(1)}\otimes
                 x^{i+1}_{(2)}\otimes\cdots\otimes\widehat{x^j}\otimes\cdots)
                 \otimes x^j
\end{align*}
By shifting the indices on the last sum and taking sum over all
$i=0,\ldots, n-1$, we get
\begin{align*}
h_{n+1}d^{CB}_n(x^1\otimes\cdots\otimes x^n)
 = & \sum_{i=1}^{n-1}\sum_{j=0}^{i-1}(-1)^{i+j+1}
               (\cdots\otimes\widehat{x^{j+1}}\otimes
                \cdots\otimes x^{i+1}_{(1)}\otimes x^{i+1}_{(2)}\otimes
                \cdots)\otimes x^{j+1}\\
   & + \rho_n(x^1\otimes\cdots\otimes x^n)\\
   & + \sum_{i=0}^{n-2}\sum_{j=i+1}^{n-1}(-1)^{i+j}
               (\cdots\otimes x^{i+1}_{(1)}\otimes
                 x^{i+1}_{(2)}\otimes\cdots\otimes\widehat{x^{j+1}}\otimes\cdots)
                \otimes x^{j+1}
\end{align*}
On the other hand
\begin{align*}
(-1)^i(\partial_i\otimes id_1)h_n(x^1\otimes\cdots\otimes x^n)
 = & \sum_{j=0}^i(-1)^{i+j+1}(\cdots\otimes\widehat{x^{j+1}}\otimes\cdots\otimes
                 x^{i+2}_{(1)}\otimes x^{i+2}_{(2)}\otimes\cdots)\otimes x^{j+1}\\
   & + \sum_{j=i+1}^{n-1}(-1)^{i+j+1}
                (\cdots\otimes x^{i+1}_{(1)}\otimes x^{i+1}_{(2)}
                 \otimes\widehat{x^{j+1}}\otimes\cdots)\otimes x^{j+1}
\end{align*}
which means
\begin{align*}
(d^{CB}_n\otimes id_1)h_n(x^1\otimes\cdots\otimes x^n)
 = & \sum_{i=1}^{n}\sum_{j=0}^{i-1}(-1)^{i+j}
                (\cdots\otimes\widehat{x^{j+1}}\otimes
                 \cdots\otimes x^{i+1}_{(1)}\otimes x^{i+1}_{(2)}
                 \otimes\cdots)\otimes x^{j+1}\\
   & + \sum_{i=0}^{n-2}\sum_{j=i+1}^{n-1}(-1)^{i+j+1}(\cdots\otimes x^{i+1}_{(1)}
                 \otimes x^{i+1}_{(2)}\otimes
                 \cdots\widehat{x^{j+1}}\otimes\cdots)\otimes x^{j+1}
\end{align*}
Then one can easily see that
\begin{align}
h_{n+1}d^{CB}_n + (d^{CB}_n\otimes id_1)h_n = \rho_n
\end{align}
for any $n\geq 1$.
\end{proof}

\section{Hochschild complex}\label{Hochschild}

\begin{rem}
Define an action of $C_n:=\left<\tau_n|\tau_n^n\right>$ the cyclic group
  of order $n$ on $C^{\otimes n}$ by letting
\begin{align}
\tau_n^{-1}(c^1\otimes\cdots\otimes c^n) 
 = (c^2\otimes\cdots\otimes c^n\otimes c^1)
\end{align}
Then observe that $\partial_j = \tau_{n+1}^j\partial_0\tau_n^{-j}$ for
any $0\leq j\leq n-1$.
\end{rem}

\begin{thm}\label{HochschildComplex}
Let $CH_*(C)=\{C^{\otimes n}\}_{n\geq 0}$ and define
\begin{align}
d^{CH}_n = d^{CB}_n + (-1)^n\tau_{n+1}^{-1}\partial_0
 = \sum_{j=0}^n(-1)^j\tau_{n+1}^j\partial_0\tau_n^{-j}
\end{align}
Then $CH_*(C)$ is a differential graded $Lie(C)$--comodule.
\end{thm}

\begin{proof}
Define $\partial_n = \tau_{n+1}^n\partial_0\tau_n^{-n}=
\tau_{n+1}^{-1}\partial_0$ and observe that for $i<j$ one still has
\begin{align}
\partial_j\partial_i(c^1\otimes\cdots\otimes c^n)
 = & \begin{cases}
     (\cdots\otimes c^{i+1}_{(1)}\otimes c^{i+1}_{(2)}\otimes
      \cdots\otimes c^j_{(1)}\otimes c^j_{(2)}\otimes
      \cdots)         & \text{ if } i<j+1<n+2\\
     (\cdots\otimes c^{i+1}_{(1)}\otimes c^{i+1}_{(2)}\otimes
      c^{i+1}_{(3)}\otimes\cdots)
                      & \text{ if } i+1=j<n+1\\
     (c^1_{(2)}\otimes\cdots\otimes c^{i+1}_{(1)}\otimes
      c^{i+1}_{(2)}\otimes\cdots\otimes c^1_{(1)})
                      & \text{ if $i<n$ and $j=n+1$}\\
     (c^1_{(3)}\otimes\cdots\otimes c^1_{(1)}\otimes c^1_{(2)})
                      & \text{ if $i=n$ and $j=n+1$}
     \end{cases}\\
 = & \partial_i\partial_{j-1}(c^1\otimes\cdots\otimes c^n)
\end{align}
Then
\begin{align*}
d^{CH}_{n+1}d^{CH}_n 
 = & \sum_{j=0}^{n+1}\sum_{i=0}^n(-1)^{i+j}\partial_j\partial_i\\
 = & \sum_{i=0}^{n}\sum_{j=0}^i(-1)^{i+j}\partial_j\partial_i
     +\sum_{i=0}^{n}\sum_{j=i+1}^{n+1}(-1)^{i+j}\partial_j\partial_i\\
 = & \sum_{i=0}^{n-1}\sum_{j=0}^i(-1)^{i+j}\partial_j\partial_i
     -\sum_{j=0}^{n-1}\sum_{i=0}^j(-1)^{i+j}\partial_i\partial_j
 =\ \  0
\end{align*}
as we wanted to show.  Recall from Theorem~\ref{BarComplex} that
$\rho_n\partial_0 = (\partial_0\otimes id_1)\rho_n$ for any $n\geq 1$.
One can also see that 
\begin{align}
\rho_n\tau_n = (\tau_n\otimes id_1)\rho_n
\end{align}
for any $n\geq 1$ too.  Then $\rho_{n+1}d^{CH}_n = (d^{CH}_n\otimes
id_1)\rho_n$ easily follows.
\end{proof}

\section{Chevalley--Eilenberg--Leibniz complex}\label{CEL}

\begin{defn}
Let $L$ be a Leibniz coalgebra.  Define a graded $k$--module
$CE_*(L)=\{L^{\otimes n}\}_{n\geq 0}$ and define two degree
$+1$ graded morphism
\begin{align}
\rho_n(l^1\otimes\cdots\otimes l^n)
  = &\sum_{i=1}^n\left(\cdots\otimes l^i_{[1]}\otimes\cdots\right)\otimes l^i_{[2]}
    & 
d^{CE}_n = & \sum_{j=1}^n(-1)^{j-1}(\rho_j\otimes id_{n-j})
\end{align}
for all $n\geq 0$.
\end{defn}

\begin{thm}\label{ChevalleyEilenbergComplex}
Define $\partial_j = (\rho_{j+1}\otimes id_{n-j-1})$ for any $0\leq
j\leq n-1$.  Then $\partial_j\partial_i=\partial_i\partial_{j-1}$ for
all $0\leq i<j$.  Therefore $d^{CE}_n =
\sum_{j=0}^{n-1}(-1)^{j}\partial_j$ is a differential on $CE_*(L)$.
\end{thm}

\begin{proof}
For $0<i<j$ and $(x^1\otimes\cdots\otimes x^n)$ from $CE_n(L)$ consider
\begin{align}
(\rho_{j+1}\otimes id_{n-j}) & (\rho_i\otimes id_{n-i})
   (x^1\otimes\cdots\otimes x^n)\nonumber\\
 = & \sum_{a=1}^i\sum_{b=1}^{a-1}(\cdots\otimes x^b_{[1]}\otimes\cdots
         \otimes x^a_{[1]}\otimes\cdots\otimes x^i\otimes
         x^a_{[2]}\otimes\cdots\otimes x^j\otimes
         x^b_{[2]}\otimes\cdots)\nonumber\\
  & + \sum_{a=1}^i(\cdots\otimes x^a_{[1][1]}\otimes\cdots\otimes x^i\otimes
         x^a_{[2]}\otimes\cdots\otimes x^j\otimes
         x^a_{[1][2]}\otimes\cdots)\label{Part1}\\
  & + \sum_{a=1}^i\sum_{b=a+1}^{i}(\cdots\otimes x^a_{[1]}\otimes\cdots
         \otimes x^b_{[1]}\otimes\cdots\otimes x^i\otimes
         x^a_{[2]}\otimes\cdots\otimes x^j\otimes
         x^b_{[2]}\otimes\cdots)\nonumber\\
  & + \sum_{a=1}^i(\cdots\otimes x^a_{[1]}\otimes\cdots\otimes x^i\otimes
         x^a_{[2][1]}\otimes\cdots\otimes x^j\otimes
         x^a_{[2][2]}\otimes\cdots)\label{Part2}\\
  & + \sum_{a=1}^i\sum_{b=i+1}^{j+1}(\cdots\otimes x^a_{[1]}\otimes\cdots
         \otimes x^i\otimes x^a_{[2]}\otimes\cdots\otimes x^{b-1}_{[1]}
         \otimes\cdots\otimes x^j\otimes x^{b-1}_{[2]}\otimes\cdots)\nonumber
\end{align}
But $L$ is a Leibniz coalgebra.  Hence 
\begin{align*}
 (x_{[1][1]}\otimes x_{[2]}\otimes x_{[1][2]}) 
+ (x_{[1]}\otimes x_{[2][1]}\otimes x_{[2][2]})
= (x_{[1][1]}\otimes x_{[1][2]}\otimes x_{[2]}) 
\end{align*}
for any $x$ in $L$.  Therefore adding (\ref{Part1}) and (\ref{Part2})
one gets
\begin{align*}
(\rho_{j+1}\otimes id_{n-j}) & (\rho_i\otimes id_{n-i})
   (x^1\otimes\cdots\otimes x^n)\\
 = & \sum_{a=1}^i\sum_{b=1}^{a-1}(\cdots\otimes x^b_{[1]}\otimes\cdots
         \otimes x^a_{[1]}\otimes\cdots\otimes x^i\otimes
         x^a_{[2]}\otimes\cdots\otimes x^j\otimes
         x^b_{[2]}\otimes\cdots)\\
  & + \sum_{a=1}^i(\cdots\otimes x^a_{[1][1]}\otimes\cdots\otimes x^i\otimes
         x^a_{[1][2]}\otimes\cdots\otimes x^j\otimes
         x^a_{[2]}\otimes\cdots)\\
  & + \sum_{a=1}^i\sum_{b=a+1}^{i}(\cdots\otimes x^a_{[1]}\otimes\cdots
         \otimes x^b_{[1]}\otimes\cdots\otimes x^i\otimes
         x^a_{[2]}\otimes\cdots\otimes x^j\otimes
         x^b_{[2]}\otimes\cdots)\\
  & + \sum_{a=1}^i\sum_{b=i+1}^j(\cdots\otimes x^a_{[1]}\otimes\cdots
         \otimes x^i\otimes x^a_{[2]}\otimes\cdots\otimes x^b_{[1]}
         \otimes\cdots\otimes x^j\otimes x^b_{[2]}\otimes\cdots)\\
 = & (\rho_i\otimes id_{n+1-i})(\rho_j\otimes id_{n-j})
   (x^1\otimes\cdots\otimes x^n)
\end{align*}
which proves $\partial_j\partial_{i-1}=\partial_{i-1}\partial_{j-1}$ for
all $0<i<j$, This is equivalent to
$\partial_j\partial_i=\partial_i\partial_{j-1}$ for all $0\leq i<j$ as
we wanted.  The proof that $d^{CE}_*$ is a differential is similar to
the proof of Theorem~\ref{HochschildComplex} where we showed $d^{CH}_*$
is a differential.
\end{proof}

\begin{thm}
$CE_*(L)$ is a differential graded $L$--comodule.  Moreover, the
$L$--coaction is null-homotopic.
\end{thm}

\begin{proof}
The comodule structure is given by the diagonal coaction $\rho_*$.  We
must show that $\rho_{n+1}d^{CE}_n = (d^{CE}_n\otimes id_1)\rho_n$ for
any $n\geq 1$.  Then by using Theorem~\ref{ChevalleyEilenbergComplex} we
get 
\begin{align*}
\rho_{n+1}d^{CE}_n 
 = & \sum_{j=0}^n(-1)^{j}\partial_n\partial_j
 = \sum_{j=0}^n(-1)^{j}\partial_j\partial_{n-1}
 = (d^{CE}_n\otimes id_1)\rho_n
\end{align*}
for any $n\geq 1$, as we wanted to show.  In order to show that the
coaction is null-homotopic, we must furnish a null-homotopy
$CE_*(L)\xra{i_*}CE_*(L)[-1]\otimes L$.  Define
\begin{align}
i_n(x^1\otimes\cdots\otimes x^n)
 = (-1)^{n+1}(x^1\otimes\cdots\otimes x^{n-1})\otimes x^n
\end{align}
Then
\begin{align*}
i_{n+1}d^{CE}_n + (d^{CE}_{n-1}\otimes id_1)i_n
 = (-1)^{n+2} d^{CE}_n + (-1)^{n+1}(d^{CE}_{n-1}\otimes id_1)
 = (-1)^n\left(d^{CE}_n - (d^{CE}_{n-1}\otimes id_1)\right)
 = \rho_n
\end{align*}
for any $n\geq 1$ as we wanted to show.
\end{proof}

\begin{rem}
One can see that $CE_*(L)$ is not only a differential graded
$L$--comodule but also a pre-cosimplicial $L$--comodule which is not a
cosimplicial $L$--comodule.  We also believe that there is no cyclic,
nor a symmetric $L$--comodule structure on $CE_*(L)$ either.
\end{rem}

\begin{lem}\label{ChevalleyIdentity}
The differential $d^{CE}_*$ satisfies the identity
\begin{align}
d^{CE}_{p+q} 
 = (d^{CE}_p\otimes id_q) + (-1)^p(id_p\otimes d^{CE}_q)
   + \sum_j (-1)^{p+j-1}(id_p\otimes\tau_{j+1}^{-1}\otimes id_{q-j})
                      (\rho_p\otimes id_{q+1})
\end{align}
for any $p,q\geq 1$.
\end{lem}

\begin{proof}
Fix $p,q\geq 1$ and consider
\begin{align*}
\rho_{p+q}(x^1\otimes\cdots\otimes x^{p+q})
 = & \sum_{j=1}^{p+q}(x^1\otimes\cdots\otimes x^j_{[1]}\otimes\cdots\otimes
                  x^{p+q})\otimes x^j_{[2]}\\
 = & \sum_{j=1}^p(id_p\otimes\tau_{q+1}^{-1})(x^1\otimes\cdots\otimes
                  x^j_{[1]}\otimes\cdots\otimes x^p\otimes x^j_{[2]}
                  \otimes x^{i+1}\otimes\cdots\otimes x^{p+q})\\
   & + \sum_{j=1}^q(x^1\otimes\cdots\otimes x^p\otimes\cdots\otimes
                      x^{p+j}_{[1]}\otimes\cdots\otimes x^{p+q})
                  \otimes x^{p+j}_{[2]}\\
 = & (id_p\otimes\tau_{q+1}^{-1})(\rho_p\otimes id_q)
                  (x^1\otimes\cdots\otimes x^{p+q})
     + (id_p\otimes\rho_q)(x^1\otimes\cdots\otimes x^{p+q})
\end{align*}
for any $x^1\otimes\cdots\otimes x^n$ from $CE_*(L)$.  Then
\begin{align*}
d^{CE}_{p+q} 
 = & \sum_{j=1}^{p+q}(-1)^{j-1}(\rho_j\otimes id_{p+q-j})\\
 = & \sum_{j=1}^p(-1)^{j-1}(\rho_j\otimes id_{p-j})\otimes id_q
      + \sum_{j=1}^q(-1)^{p+j-1}(\rho_{p+j}\otimes id_{q-j})\\
 = & (d^{CE}_p\otimes id_q) + (-1)^p(id_p\otimes d^{CE}_q)
      + \sum_j (-1)^{p+j-1}(id_p\otimes\tau_{j+1}^{-1}\otimes id_{q-j})
              (\rho_p\otimes id_{q+1})
\end{align*}
as we wanted to show.
\end{proof}

\begin{defn}
Define a new differential graded $L$--comodule by letting
$CE^\text{red}_*(L):= CE_*(L)^L = ker(\rho_*)$.  We will call this
complex as the reductive Chevalley--Eilenberg--Leibniz complex.
\end{defn}

\begin{thm}
$CE^\text{red}_*(L)$ is a differential graded algebra with tensor
multiplication being the underlying product structure.
\end{thm}

\begin{proof}
Take $\Psi$ from $CE^\text{red}_p(L)$ and $\Phi$ from
$CE^\text{red}_q(L)$ arbitrary. By Lemma~\ref{ChevalleyIdentity} we see
that
\begin{align}
\rho_{p+q}(\Psi\otimes\Phi)
 = & (id_p\otimes\tau_{q+1}^{-1})(\rho_p(\Psi)\otimes\Phi)
     + (\Psi\otimes\rho_q(\Phi))
\end{align}
which implies $\Psi\otimes\Phi$ is in $CE^\text{red}_{p+q}(L)$. The fact
that $\rho_p(\Psi)=0$ implies
\begin{align}
d^{CE}_{p+q}(\Psi\otimes\Phi)
 = & (d^{CE}_p(\Psi)\otimes\Phi) + (-1)^p(\Psi\otimes d^{CE}_q(\Phi))
     + \sum_j (-1)^{p+j-1}(id_p\otimes\tau_{j+1}^{-1}\otimes id_{q-j})
              (\rho_p(\Psi)\otimes\Phi)\\
 = & (d^{CE}_p(\Psi)\otimes\Phi) + (-1)^p(\Psi\otimes d^{CE}_q(\Phi))
\end{align}
as we wanted to show.
\end{proof}

\section{Tale of three complexes}\label{Tale}

\begin{lem}\label{CBtoCH}
For any $n\geq 0$, one has $(id_{n+1}-(-1)^n\tau_{n+1}^{-1})d^{CH}_n =
d^{CB}_n(id_n-(-1)^{n-1}\tau_n^{-1})$.
\end{lem}

\begin{proof}
The proof is by direct calculation. Observe that
$\tau_{n+1}^{-1}\partial_i = \partial_{i-1}\tau_n^{-1}$ for any $0<i\leq
n$ and $\tau_{n+1}^{-1}\partial_0=\partial_n$.  Consider
\begin{align*}
\left(id_{n+1}-(-1)^n\tau_{n+1}^{-1}\right)d^{CH}_n 
 = & \sum_{i=0}^n(-1)^i\partial_i 
   - \sum_{i=0}^n(-1)^{n+i}\tau_{n+1}^{-1}\partial_i \\
 = & \sum_{i=0}^n(-1)^i\partial_i 
   - \sum_{i=1}^n(-1)^{n+i}\partial_{i-1}\tau_n^{-1} 
   - (-1)^n\partial_n\\
 = & \sum_{i=0}^{n-1}(-1)^i\partial_i 
   - \sum_{i=0}^{n-1}(-1)^{n-1+i}\partial_i\tau_n^{-1}\\
 = & d^{CB}_n\left(id_n-(-1)^{n-1}\tau_n^{-1}\right)
\end{align*}
which proves what we wanted.
\end{proof}

\begin{lem}\label{CHtoCB}
Let $N_n=\sum_{j=0}^{n-1}(-1)^{(n-1)j}\tau_n^j$ for any $n\geq 1$.  Then
$N_{n+1}d^{CB}_n = d^{CH}_nN_n$.
\end{lem}

\begin{proof}
Proof is by direct calculation.  Since $\tau_{n+1}^{-n-1}
\partial_0\tau_n^n =\partial_0$, we have
\begin{align}
\partial_i\tau_n^{-j}
 = & \begin{cases}
     \tau_{n+1}^{-j}  \partial_{i+j}   & \text{ if } i+j<n\\
     \tau_{n+1}^{-j-1}\partial_{i+j-n} & \text{ if } i+j\geq n
     \end{cases}
\end{align}
Then
\begin{align*}
d^{CH}_nN_n
 = & \sum_{j=0}^{n-1}\sum_{i=0}^n(-1)^{(n-1)j+i}\partial_i\tau_n^{-j}\\
 = & \sum_{j=0}^{n-1}\sum_{i=0}^{n-j-1}(-1)^{(n-1)j+i}\tau_{n+1}^{-j}\partial_{i+j}
  +  \sum_{j=0}^{n-1}\sum_{i=n-j}^n(-1)^{(n-1)j+i}\tau_{n+1}^{-j-1}\partial_{i+j-n}\\
 = & \sum_{j=0}^{n-1}\sum_{i=j}^{n-1}(-1)^{(n-1)j+i+j}\tau_{n+1}^{-j}\partial_i
  + \sum_{j=0}^{n-1}\sum_{i=0}^j(-1)^{(n-1)j+i+j-n}\tau_{n+1}^{-j-1}\partial_i
\end{align*}
By switching the order of summation we get
\begin{align*}
d^{CH}_nN_n
 = & \sum_{i=0}^{n-1}\sum_{j=0}^{i}(-1)^{nj+i}\tau_{n+1}^{-j}\partial_i
  +  \sum_{i=0}^{n-1}\sum_{j=i}^{n-1}(-1)^{n(j+1)+i}\tau_{n+1}^{-j-1}\partial_i\\
 = & \sum_{i=0}^{n-1}\sum_{j=0}^{i}(-1)^{nj+i}\tau_{n+1}^{-j}\partial_i
  + \sum_{i=0}^{n-1}\sum_{j=i+1}^{n}(-1)^{nj+i}\tau_{n+1}^{-j}\partial_i\\
 = & N_{n+1}d^{CB}_n
\end{align*}
as we wanted to show.
\end{proof}

\begin{thm}\label{CyclicComplex}
Let $t_* = (id_* - (-1)^*\tau_*^{-1})$.  Then there is an exact periodic
sequence of differential graded $Lie(C)$--comodules of the form
\begin{align}
CH_*(C)\xra{t_*}CB_*(C)\xra{N_*}CH_*(C)\xra{t_*}CB_*(C)\xra{N_*}\cdots
\end{align}
\end{thm}

\begin{proof}
Lemma~\ref{CBtoCH} shows that $CB_*(C)\xra{t_*}CH_*(C)$ is a morphism of
differential graded $k$--modules and Lemma~\ref{CHtoCB} shows that
$CH_*(C)\xra{N_*}CB_*(C)$ is a morphism differential graded
$k$--modules.  One can easily see that $t_*N_*=N_*t_*=0$.  So, in order
to prove exactness, we must show that $ker(t_*)=im(N_*)$ and
$ker(N_*)=im(t_*)$.  Assume $\bf{c}$ is in $ker(t_n)$ which means
$\tau_n{\bf c}=(-1)^{n-1}{\bf c}$.  Then
\begin{align*}
N_n{\bf c} 
= \sum_{j=0}^{n-1}(-1)^{(n-1)j}\tau_n^j{\bf c}
= \sum_{j=0}^{n-1}(-1)^{(n-1)j}(-1)^{(n-1)j}{\bf c}
= n{\bf c}
\end{align*}
Hence ${\bf c}$ is in $im(N_*)$ since we assume $k$ is a field of
characteristic $0$. Conversely, assume ${\bf c}$ is in $ker(N_*)$.  Then 
\begin{align*}
n{\bf c} 
  = n{\bf c} - N_*{\bf c}
  = \sum_{j=1}^{n-1}\left(id_n-(-1)^{(n-1)j}\tau_n^j\right){\bf c}
  = \left(id_n-(-1)^{(n-1)}\tau_n\right)
    \sum_{j=1}^{n-1}\sum_{i=0}^j(-1)^{(n-1)i}\tau_n^i{\bf c}
\end{align*}
which means ${\bf c}$ is in the image of $t_*$.  The fact that both
$t_*$ and $N_*$ are morphism of differential graded $Lie(C)$--comodules
follows from the Lemma~\ref{RepresentationCommute}.
\end{proof}

\begin{rem}
The action of $C_n$ on $C^{\otimes n}$ also extends to defining an
action of $\Sigma_n$ on $C^{\otimes n}$.  All we need to define is the
action of a single transposition.  For $n\geq 3$ define
\begin{align}
(\tau_2\otimes id_{n-2})\cdot(x^1\otimes x^2\otimes\cdots\otimes x^n)
 = (x^2\otimes x^1\otimes\cdots\otimes x^n)
\end{align}
for any $(x^1\otimes\cdots\otimes x^n)$ from $C^{\otimes n}$ and $n\geq
3$.
\end{rem}

\begin{lem}\label{RepresentationCommute}
Given any $\sigma\in\Sigma_n$ for any $n\geq 1$, one has $\rho_n\sigma =
(\sigma\otimes id_1)\rho_n$.
\end{lem}

\begin{lem}\label{Identity}
Define an element $h_n =
\sum_{j=0}^{n-1}(-1)^{j+1}(id_j\otimes\tau_{n-j}^{-1})$ in $k[\Sigma_n]$
for any $n\geq 1$.  Define also
$\epsilon_n=\sum_{\sigma\in\Sigma_n}\sgn(\sigma)\sigma$ in
$k[\Sigma_n]$.  Then one has
\begin{align}
(\epsilon_n\otimes id_1)h_{n+1} = (-1)^{n+1}\epsilon_{n+1}
\end{align}
for any $n\geq 1$.
\end{lem}

\begin{proof}
One can consider $\Sigma_n\subseteq \Sigma_{n+1}$ as the stabilizer of
the element $n$ if one considers the symmetric group on $n+1$ letters as
the group of automorphism of the finite set $\{0,\ldots,n\}$.  Then we
claim
\begin{align}
\Sigma_{n+1} = \bigsqcup_{j=0}^n\Sigma_n\tau_{n+1}^{j}
\end{align}
If we assume on the contrary that
$\Sigma_n\tau_{n+1}^{i}=\Sigma_n\tau_{n+1}^{j}$ for $i\neq j$, then one
necessarily has $\tau_{n+1}^{i-j}\in\Sigma_n$ which is not the case
since $\tau_{n+1}^{i-j}(n)=n$ only if $i=j$.  Now observe that
\begin{align}
\tau_{n+1}^{-j} = (\tau_n^{-j+1}\otimes id_1)(id_j\otimes\tau_{n+1-j}^{-1})
\end{align}
for any $1\leq j\leq n$.  This means
\begin{align}
\Sigma_{n+1} 
= \bigsqcup_{j=0}^n\Sigma_n(\tau_n^{-j+1}\otimes id_1)(id_j\otimes\tau_{n+1-j}^{-1})
= \bigsqcup_{j=0}^n\Sigma_n(id_j\otimes\tau_{n+1-j}^{-1})
\end{align}
which finishes the proof.
\end{proof}

\begin{thm}\label{AntiCommuteMorphism}
There is a morphism of differential graded $Lie(C)$--modules of the form
$CB_*(C)\xra{\epsilon_*}CE_*(Lie(C))$.
\end{thm}

\begin{proof}
We must show that $\epsilon_{n+1}d^{CB}_n=d^{CE}_n\epsilon_n$ for any
$n\geq 1$.  The proof is by induction on $n$.  First observe that
$\partial_0^{CE}(x) = \delta(x) = x_{(1)}\otimes x_{(2)} -
x_{(2)}\otimes x_{(1)}$ by the definition of $Lie(C)$ and the cobracket
associated with $Lie(C)$.  Then for $n=1$ we have
\begin{align*}
\epsilon_2d^{CB}_1 
 = & (id_2-\tau_2)\partial^{CB}_0 = \delta_1 = d^{CE}_1\epsilon_1
\end{align*}
trivially satisfied.  Assume as the induction hypothesis that
$\epsilon_{n+1}d^{CB}_n = d^{CE}_n\epsilon_n$.  By using
Lemma~\ref{Identity} and Theorem~\ref{BarHomotopy} we get
\begin{align*}
\epsilon_{n+2}d^{CB}_{n+1}
 = & (-1)^{n+2}(\epsilon_{n+1}\otimes id_1)h_{n+2}d^{CB}_{n+1}\\
 = & (-1)^{n+2}(\epsilon_{n+1}\otimes id_1)
        \left(-(d^{CB}_n\otimes id_1)h_{n+1} 
              +\rho_{n+1}\right)\\
 = & (-1)^{n+1}(d^{CE}_n\otimes id_1)(\epsilon_n\otimes id_1)h_{n+1}
        + (-1)^{n}(\epsilon_{n+1}\otimes id_1)\rho_{n+1}
\end{align*}
Lemma~\ref{RepresentationCommute} implies $(\epsilon_{n+1}\otimes
id_1)\rho_{n+1} = \rho_{n+1}\epsilon_{n+1}$.  Then
\begin{align*}
\epsilon_{n+2}d^{CB}_{n+1}
 = & (d^{CE}_n\otimes id_1)\epsilon_{n+1} 
     + (-1)^{n}\rho_{n+1}\epsilon_{n+1}
 = d^{CE}_{n+1}\epsilon_{n+1}
\end{align*}
for any $n\geq 1$, as we wanted to show. This proves there is a morphism
of differential graded $Lie(C)$--comodules of the form
$CB_*(C)\xra{\epsilon_*}CE_*(Lie(C))$.  The fact that $\epsilon_*$ is a
morphism of $Lie(C)$--comodules follows from
Lemma~\ref{RepresentationCommute}.
\end{proof}

\begin{thm}\label{CHtoCE}
There is a morphism of differential graded $Lie(C)$--comodules of the
form $CH_*(C)\xra{id_1\otimes\epsilon_{*-1}}CE_*(Lie(C))$.
\end{thm}

\begin{proof}
We must show that $d^{CH}_n(id_1\otimes\epsilon_{n-1})=
(id_1\otimes\epsilon_n)d^{CE}_n$ for $n\geq 1$.  For $n=1$, one has
\begin{align*}
d^{CH}_1(x) 
= (x_{(1)}\otimes x_{(2)}) - (x_{(2)}\otimes x_{(1)}) 
= \rho_1(x)
= d^{CE}_1(x)
\end{align*}
For $n\geq 2$ consider
\begin{align*}
(id_1\otimes\epsilon_{n+1})d^{CH}_{n+1}
 = & (id_1\otimes\epsilon_{n+1})\partial_0
   - (id_1\otimes\epsilon_{n+1})(id_1\otimes d^{CB}_n)
   + (-1)^{n+1}(id_1\otimes\epsilon_{n+1})\tau_{n+2}^{-1}\partial_0\\
 = & - (id_1\otimes d^{CE}_n)(id_1\otimes\epsilon_n)
     + (id_1\otimes\epsilon_{n+1})\partial_0
     - (-1)^n(id_1\otimes\epsilon_{n+1})(id_1\otimes\tau_{n+1}^{-1})
            (\tau_2\otimes id_n)\partial_0\\
 = & - (id_1\otimes d^{CE}_n)(id_1\otimes\epsilon_n)
     + (id_1\otimes\epsilon_{n+1})\partial_0
     - (id_1\otimes\epsilon_{n+1})(\tau_2\otimes id_n)\partial_0\\
 = & - (id_1\otimes d^{CE}_n)(id_1\otimes\epsilon_n)
     + (id_1\otimes\epsilon_{n+1})(\rho_1\otimes id_n)
\end{align*}
by using the fact that $d^{CE}_n\epsilon_n=\epsilon_{n+1}d^{CB}_n$ for
any $n\geq 0$.  Then by using Lemma~\ref{ChevalleyIdentity} and
Lemma~\ref{Identity} we get
\begin{align*}
(id_1\otimes\epsilon_{n+1})d^{CH}_{n+1}
 = & - (id_1\otimes d^{CE}_n)(id_1\otimes\epsilon_n) 
     + (id_2\otimes\epsilon_n)(d^{CE}_1\otimes id_n)\\ 
   & + \sum_{j=1}^n(-1)^{j}(id_1\otimes\tau_{j+1}^{-1}\otimes id_{n-j})
       (id_2\otimes\epsilon_n)(\rho_1\otimes id_n)\\ 
 = &\ (d^{CE}_1\otimes id_n)(id_1\otimes\epsilon_n) 
     - (id_1\otimes d^{CE}_n)(id_1\otimes\epsilon_n)\\ 
   & + \sum_{j=1}^n(-1)^{j}(id_1\otimes\tau_{j+1}^{-1}\otimes id_{n-j})
       (\rho_1\otimes id_n)(id_1\otimes\epsilon_n)\\ 
 = &\ d^{CE}_{n+1}(id_1\otimes\epsilon_n)
\end{align*}
by using the fact that 
\begin{align}
\epsilon_{n+1} = & \sum_{j=0}^n
(-1)^j(\tau_{j+1}^{-1}\otimes id_{n-j})(id_1\otimes\epsilon_n)
\end{align}
whose proof is similar to Lemma~\ref{Identity}.  The fact that the
morphism we defined is a morphism of differential graded
$Lie(C)$--comodules follows from Lemma~\ref{RepresentationCommute}.
\end{proof}

\begin{cor}\label{AllTogether}
There is a commutative diagram of differential graded
$Lie(C)$--comodules of the form
\begin{equation}
\begin{CD}
CB_*(C)                    @ =                  CB_*(C)\\
@V{N_*}VV                                       @VV{\epsilon_*}V\\
CH_*(C)       @>>{(id_1\otimes\epsilon_{*-1})}> CE_*(Lie(C))
\end{CD}
\end{equation}
\end{cor}
The algebra version of Corollary~\ref{AllTogether} is proven by by
J.M. Lodder in \cite{Lodder:FromLeibnizToCyclicHomology}.

\begin{lem}\label{Ideal}
$ker(\epsilon_*)$ is a two sided ideal of $CB_*(C)$.
\end{lem}

\begin{proof}
Take $\Psi$ from $ker(\epsilon_n)$ and $\Phi$ from $CB_m(C)$ arbitrary.
Consider
\begin{align}
n!\epsilon_{n+m}(\Psi\otimes\Phi)
 = \epsilon_{n+m}(\epsilon_n\otimes id_m)(\Psi\otimes\Phi)
 = 0
\end{align}
So, $\Psi\otimes\Phi$ belongs to $ker(\epsilon_*)$,
i.e. $ker(\epsilon_*)$ is a right ideal.  The proof that
$ker(\epsilon_*)$ is a left ideal is similar.
\end{proof}

\begin{defn}
Define a graded commutator by letting
\begin{align}
[\Psi,\Phi] = (\Psi\otimes\Phi) + (-1)^{nm}(\Phi\otimes\Psi)
\end{align}
for any $\Psi\in CB_n(C)$ and $\Phi\in CB_m(C)$.  Let
$[CB_*(C),CB_*(C)]$ denote the graded \underline{sub-module} of graded
commutators. 
\end{defn}

\begin{lem}
$[CB_*(C),CB_*(C)]=ker(N_*)$.  Therefore $[CB_*(C),CB_*(C)]$ is a
differential graded $Lie(C)$--sub-comodule of $CB_*(C)$.
\end{lem}

\begin{proof}
For any $\Psi\in CB_n(C)$ and $\Phi\in CB_m(C)$ we have
\begin{align*}
(-1)^{(n+m-1)n}\tau_{n+m}^{-n}(\Psi\otimes\Phi)
= & (-1)^{nm}(\Phi\otimes\Psi)
\end{align*}
which means $[\Psi,\Phi]$ is in the image of
$(id_{n+m}-(-1)^{(n+m-1)n}\tau_{n+m}^{-n})$ which is in $ker(N_*)$ since
\begin{align}
(-1)^{n+m-1}\tau_{n+m}N_{n+m}  = N_{n+m}
\end{align}
for any $n,m$ from $\B{N}$.  Conversely, assume $\Psi$ is in $ker(N_*)$
which is equal to $im(t_*)$ by Theorem~\ref{CyclicComplex}.  Then $\Psi$
is of the form $\Psi = (id_n - (-1)^{n-1}\tau_n)\Psi'$ for some $\Psi'$
from $CB_*(C)$.  Anything which lies in the image of
$(id_n-(-1)^{n-1}\tau_n)$ lies in the submodule $[CB_*(C),CB_*(C)]$.
This finishes the proof.
\end{proof}

\begin{lem}
$[CB_*(C),CB_*(C)]$ is contained in $ker(\epsilon_*)$.
\end{lem}

\begin{proof}
The ideal $ker(\epsilon_*)$ is generated by elements of the form
$(1-\sgn(\sigma)\sigma)\Psi$ where $\Psi$ is from $CB_n(C)$ and $\sigma$
is from $\Sigma_n$.  Moreover $ker(N_*)=im(t_*)$ thus is generated by
elements of the form $(1-\sgn(\tau_n)\tau_n)$.
\end{proof}

\begin{defn}
Define a new differential graded $Lie(C)$--comodule $CC^\lambda_*(C)$ by
letting
\begin{align*}
CC^\lambda_*(C)[+1] := CB_*(C)/[CB_*(C),CB_*(C)] \cong im(N_*)
\end{align*}
We will call this complex as the cyclic complex associated with the
coassociative coalgebra $C$.
\end{defn}

\begin{rem}
A coassociative coalgebra is called $H$--counital if $CB_*(C)$ has
trivial homology.  As one can observe, if $C$ is counital then it is
$H$--counital.  Moreover, If $C$ is $H$--counital, then the double
complex we described above in Theorem~\ref{CyclicComplex} is homotopic
to $CC^\lambda_*(C)$,
\end{rem}

\begin{defn}
Define a new differential graded $Lie(C)$--comodule
$CE^\text{sym}_*(Lie(C))$ by letting
\begin{align*}
CE^\text{sym}_*(Lie(C)):= CB_*(C)/ker(\epsilon_*) \cong im(\epsilon_*)
\end{align*}
We will call this complex as the Chevalley--Eilenberg--Lie complex
associated with the Lie coalgebra $Lie(C)$.
\end{defn}

\begin{defn}
Let $\Sigma_{(p,q)}$ be the set of $(p,q)$--shuffles in $\Sigma_{p+q}$
\cite[Appendix A]{Loday:CyclicHomology}.  Define two elements from
$k[\Sigma_{p+q}]$ as
\begin{align}
\epsilon^{(p,q)} = & \sum_{\sigma\in\Sigma_{(p,q)}}\sgn(\sigma)\sigma\\
\epsilon_{(p,q)} = & \sum_{\sigma\in\Sigma_{(p,q)}}\sgn(\sigma)\sigma^{-1}
\end{align}
for any $p,q\geq 1$.
\end{defn}

\begin{lem}
For any $p,q\geq 1$ one has
$\epsilon_{p+q}=(\epsilon_p\otimes\epsilon_q)\epsilon^{(p,q)}$ and
therefore $\epsilon_{p+q}=\epsilon_{(p,q)}(\epsilon_p\otimes\epsilon_q)$ 
\end{lem}

\begin{proof}
We claim that there is a unique element $\epsilon^{(p,q)}\in
k[\Sigma_{p+q}]$ determined by the following two conditions: (i)
$(\epsilon_p\otimes\epsilon_q)\epsilon^{(p,q)}=\epsilon_{p+q}$ and more
importantly (ii) $\sigma(1)<\cdots<\sigma(p)$ and
$\sigma(p+1)<\cdots<\sigma(p+q)$ for any $\sigma$ appearing as a summand
in $\epsilon^{(p,q)}$.  Consider $\Sigma_p\times\Sigma_q$ as a subgroup
of $\Sigma_{p+q}$ and consider its right coset space
$(\Sigma_p\times\Sigma_q)\backslash\Sigma_{p+q}$.  The element
$\epsilon^{(p,q)}$ will be the sum $\sum_\sigma \sgn(\sigma)\sigma$
where the sum is taken over unique representatives $\sigma$ from each
right coset of $\Sigma_p\times\Sigma_q$.  So the first condition is
satisfied.  Assume $(\Sigma_p\times\Sigma_q)\sigma$ be a right coset and
consider the following disjoint union 
\begin{equation*}
\{\sigma(1),\ldots,\sigma(p)\}\sqcup\{\sigma(p+1),\ldots,\sigma(p+q)\}
\end{equation*}
of the set $\{1,\ldots, n\}$.  One must see that, there is a unique
element $(\gamma,\gamma')\in\Sigma_p\times\Sigma_q\subseteq\Sigma_{p+q}$
such that
\begin{align*}
\sigma(\gamma(1)) & < \cdots < \sigma(\gamma(p))\\
\sigma(\gamma'(1)+p)& < \cdots < \sigma(\gamma'(q)+p)
\end{align*}
which proves that every element $\sigma$ can be chosen uniquely from
$\Sigma_{(p,q)}$ thus
\begin{align}
(\Sigma_p\times\Sigma_q)\Sigma_{(p,q)} = \Sigma_{p+q}
\end{align}
This finishes the proof.
\end{proof}

\begin{defn}
Let $X$ be a $k$--module.  Define the shuffle product on the graded
$k$--module $\{X^{\otimes n}\}_{n\geq 1}$ as
\begin{align}
\mu^{sh}_{p,q}(\Psi\otimes\Phi) = \epsilon_{(p,q)}(\Psi\otimes\Phi)
\end{align}
for any $\Psi$ from $X^{\otimes p}$ and $\Phi$ from $X^{\otimes q}$.
\end{defn}

\begin{thm}\label{CommutativeAlgebra}
$CE^\text{sym}_*(Lie(C))$ is a differential graded commutative algebra
with respect to the shuffle product. Moreover,
$CB_*(C)\xra{\epsilon_*}CE^\text{sym}_*(Lie(C))$ is an epimorphism of
differential graded algebras.
\end{thm}

\begin{proof}
Take $\epsilon_p(\Psi)$ from $CE^\text{sym}_p(Lie(C))$ and
$\epsilon_q(\Phi)$ from $CE^\text{sym}_q(Lie(C))$.  Consider
\begin{align*}
\epsilon_{(p,q)}(\epsilon_p(\Psi)\otimes\epsilon_q(\Phi))
 = \epsilon_{(p,q)}(\epsilon_p\otimes\epsilon_q)(\Psi\otimes\Phi)
 = \epsilon_{p+q}(\Psi\otimes\Phi)
\end{align*}
for any $\Psi$ from $CB_p(C)$ and $\Phi$ from $CB_q(C)$.  This means the
shuffle product of any two element from $CE^\text{sym}_*(Lie(C))$ is
again in $CE^\text{sym}_*(Lie(C))$.  Moreover, since we have
\begin{align*}
\epsilon_{(q,p)}(\epsilon_q(\Phi)\otimes\epsilon_p(\Psi))
 = \epsilon_{p+q}(\Phi\otimes\Psi)
 = \epsilon_{p+q}\tau_{p+q}^{-q}(\Psi\otimes\Phi)
 = (-1)^{pq}\epsilon_{(p,q)}(\epsilon_p(\Psi)\otimes\epsilon_q(\Phi))
\end{align*}
one can see that the shuffle product is graded commutative.  Now,
consider
\begin{align*}
d^{CE}_{p+q}\epsilon_{(p,q)}(\epsilon_p(\Psi)\otimes\epsilon_q(\Phi))
 = & d^{CE}_{p+q}\epsilon_{p+q}(\Psi\otimes\Phi)\\
 = & \epsilon_{p+q+1}d^{CB}_{p+q}(\Psi\otimes\Phi)\\
 = & \epsilon_{p+q+1}(d^{CB}_p(\Psi)\otimes\Phi)
    +(-1)^p\epsilon_{p+q+1}(\Psi\otimes d^{CB}_q(\Phi))\\
 = & \epsilon_{(p+1,q)}(\epsilon_{p+1}d^{CB}_p(\Psi)\otimes\epsilon_q\Phi)
    +(-1)^p\epsilon_{(p,q+1)}(\epsilon_p(\Psi)\otimes
                              \epsilon_{q+1}d^{CB}_q(\Phi))\\
 = & \epsilon_{(p+1,q)}(d^{CE}_p\epsilon_p(\Psi)\otimes\epsilon_q(\Phi))
     +(-1)^p\epsilon_{(p,q+1)}(\epsilon_p(\Psi)\otimes d^{CE}_q\epsilon_q(\Phi))
\end{align*}
as we wanted to show.  This finishes the proof.
\end{proof}

\begin{thm}\label{CocommutativeCoalgebra}
$CE^\text{sym}_*(Lie(C))$ is a graded cocommutative coalgebra with
respect to the deconcatenation coproduct.
\end{thm}

\begin{proof}
Recall that
$\epsilon_{p+q}=(\epsilon_p\otimes\epsilon_q)\epsilon^{(p,q)}$.  This
means
\begin{align*}
\Delta(\epsilon_n\Psi)
 = & \sum_{p+q=n}\Delta_{(p,q)}(\epsilon_p\otimes\epsilon_q)
                 \epsilon^{(p,q)}\Psi \\
 = & \sum_{p+q=n}\sum_{\Phi_{(1)}\otimes\Phi_{(2)}=\epsilon^{(p,q)}\Psi}
                 \epsilon_p\Phi_{(1)}\otimes\epsilon_q\Phi_{(2)}
\end{align*}
implying $\Delta(CE^\text{sym}_*(Lie(C)))\subseteq
CE^\text{sym}_*(Lie(C))\otimes CE^\text{sym}_*(Lie(C))$ as we wanted to
show.  Moreover
\begin{align*}
\Delta^{op}(\epsilon_n\Psi)
 = & \sum_{p+q=n}\sum_{\Phi_{(1)}\otimes\Phi_{(2)}=\epsilon^{(p,q)}\Psi}
     \epsilon_q\Phi_{(2)}\otimes\epsilon_p\Phi_{(1)}
\end{align*}
Since $\Phi_{(2)}\otimes\Phi_{(1)}=
\tau_{p+q}^{-p}(\Phi_{(1)}\otimes\Phi_{(2)})$ and
$\tau_{p+q}^{-p}\epsilon^{(p,q)}\tau_{p+q}^p=\epsilon^{(q,p)}$ one has
\begin{align*}
\Delta^{op}(\epsilon_n\Psi)
 = & \sum_{p+q=n}\sum_{\Phi_{(2)}\otimes\Phi_{(1)}=\epsilon^{(q,p)}
                       \tau_{p+q}^{-p}(\Psi)}
     \left(\epsilon_q\Phi_{(2)}\otimes\epsilon_p\Phi_{(1)}\right)\\
 = & \Delta(\epsilon_n\tau_{p+q}^{-p}\Psi)\\
 = & (-1)^{pq}\Delta(\epsilon_n\Psi)
\end{align*}
which proves the graded cocommutativity.
\end{proof}

\begin{rem}
According to \cite[Appendix A]{Loday:CyclicHomology} the collection
$CE^\text{sym}_*(Lie(C))$ is a graded Hopf algebra with respect to the
multiplication and comultiplication structures we defined above.
According to \cite{MilnorMoore:HopfAlgebras}, the graded Hopf algebra
$CE^\text{sym}_*(Lie(C))$ is generated by the submodule of primitive
elements as an exterior algebra.  The submodule of primitive elements
consists of degree 1 elements which fails to be a differential graded
submodule.  Thus, this result can not be extended to a generation result
as a {\bf differential} graded Hopf algebra.
\end{rem}

\section{Miscellanea}\label{Misc}

\begin{defn}
Let $M_n^c(k)$ be the free $k$--module generated by symbols of the form
$\{e_{ij}\}_{i,j=1}^n$.  The comultiplication on the generators are
defined as
\begin{align}
\Delta(e_{ij}) = \sum_a e_{ia}\otimes e_{aj}
\end{align}
for any $e_{ij}$ from the basis.  One can easily see that
\begin{align*}
(\Delta\otimes id)\Delta(e_{ij})
 = \sum_{a,b} e_{ia}\otimes e_{ab}\otimes e_{bj}
 = (id\otimes\Delta)\Delta(e_{ij})
\end{align*}
for any $e_{ij}$ from the basis, hence $M_n^c(k)$ is a coassociative
coalgebra.  Define also a morphism $M_n^c(k)\xra{\eta}k$ by
\begin{align}
\eta(e_{ij}) 
 = \delta^i_j 
 = \begin{cases}
   1  & \text{ if } i=j\\
   0  & \text{ otherwise }
   \end{cases}
\end{align}
Then 
\begin{align*}
(\eta\otimes id)\Delta(e_{ij})
 = \sum_a \delta^i_ae_{aj} = e_{ij}
 = \sum_a e_{ia}\delta^a_j
 = (id\otimes\eta)\Delta(e_{ij})
\end{align*}
for any $e_{ij}$ which means $\eta$ is a counit.
\end{defn}

\begin{thm}
Let $X$ and $Y$ be two coassociative counital coalgebras.  Then there is
a coassociative counital coalgebra structure on $X\otimes Y$.
\end{thm}

\begin{proof}
For any $x\otimes y$ from $X\otimes Y$ define
\begin{align}
\Delta(x\otimes y)
 = & (x_{(1)}\otimes y_{(1)})\otimes(x_{(2)}\otimes y_{(2)})\\
\eta(x\otimes y)
 = & \eta(x)\eta(y)
\end{align}
It is easy to check that $X\otimes Y$ is a coassociative counital coalgebra.
\end{proof}

\begin{cor}
Let $C$ be a coassociative counital coalgebra.  Then one has
coassociative counital coalgebras of the form $M_n^c(C)$ defined by
$M_n^c(k)\otimes C$ for any $n\geq 1$.
\end{cor}

\begin{cor}
Let $C$ be a coassociative coalgebra.  Then one has Lie coalgebras of
the form $gl_n^c(C)$ defined by $Lie(M_n^c(C))$ for any $n\geq 1$.
\end{cor}

\begin{defn}
Given any $k$--module $X$, define $X^\vee$ as $\text{Hom}_k(X,k)$.
\end{defn}

\begin{lem}
Let $C$ be a coassociative (Leibniz) coalgebra.  Then $C^\vee$ is a
associative (resp. Leibniz) algebra.
\end{lem}

\begin{proof}
Let $C\xra{\delta}C\otimes C$ be the underlying coassociative (Leibniz)
coalgebra structure.  Then one has
\begin{align}
\text{Hom}_k(C,k)\otimes \text{Hom}_k(C,k)\xra{\mu}
\text{Hom}_k(C\otimes C,k)\xra{\text{Hom}_k(\delta,k)}\text{Hom}_k(C,k)
\end{align}
where $\mu$ is the pointwise multiplication map.  The fact that $C^\vee$
is associative follows from the fact that $C$ is coassociative and
$\mu$ is associative.  If $C$ is a Leibniz coalgebra, then given any
three $f,g,h$ from $C^\vee$, one has
\begin{align}
([[f,g],h]-[[f,h],g])(x) 
= &   f(x_{[1][1]})g(x_{[1][2]})h(x_{[2]})
    - f(x_{[1][1]})g(x_{[2]})h(x_{[1][2]})\\
= & f(x_{[1]})g(x_{[2]]1]})h(x_{[2][2]})\\
= & [f,[g,h]](x)
\end{align}
for any $x\in C$ proving $C^\vee$ is a Leibniz algebra.
\end{proof}

\begin{thm}
Let $C$ be a finite dimensional coassociative (Leibniz) coalgebra.  Then
the category of finite dimensional $C$--comodules is anti-equivalent to
the category of finite dimensional $C^\vee$--modules through the functor
$\text{Hom}_k(\cdot,k)$.
\end{thm}

\begin{cor}
Let $L$ be a finite dimensional Leibniz coalgebra.  Then for any finite
dimensional right $L$--comodule $X$, one has $X^L\cong
(X^\vee)_{L^\vee}$.
\end{cor}

\begin{proof}
Let $X\xra{\rho_X}X\otimes L$ be the underlying coalgebra structure.
$X^L$ is defined as $X^L = ker(\rho_X)$.  For a Leibniz algebra $G$ and
a right $G$--module $Y\otimes G\xra{\rho_Y}Y$, the module $Y_G$ is
defined as $coker(\rho_Y)$.  Since $(\cdot)^\vee$ is exact, 
\begin{align}
\left(ker(\rho_X)\xra{}X\xra{\rho_X}X\otimes L\right)^\vee
 = \left(coker(\rho_{X^\vee})\xla{}X^\vee\xla{\rho_{X^\vee}}
         X^\vee\otimes L^\vee\right)
\end{align}
for any finite dimensional $L$--comodule $X$, as we wanted to show.
\end{proof}

\begin{cor}\label{DualIsomorphism}
Let $M_n^a(k)$ be the algebra of $n\times n$--matrices.  Then
$\text{Hom}_k(CB_*(M_n^a(k)),k)\cong CB_*(M_n^c(k))$.  Moreover,
$\text{Hom}_k(CB_*(M_n^a(k))_{gl_n(k)},k)\cong CB_*(M_n^c(k))^{gl_n^c(k)}$
\end{cor}

\begin{cor}\label{Coreductivity}
Any finite dimensional $gl_n^c(k)$--comodule $X$ splits as
$X^{gl_n^c(k)}\oplus X^\text{rest}$.
\end{cor}

\begin{proof}
The dual $gl_n^c(k)^\vee$ is the Lie algebra $gl_n(k)$.  And since
$gl_n(k)$ is a reductive Lie algebra, the result follows.
\end{proof}

\section{Loday--Quillen--Tsygan Theorem}\label{LQT}

The differential graded module $CB_*(C)$ is actually a cosimplicial
$k$--module.  Thus, the total complex of the bicosimplicial $k$--module
$CB_*(C)\otimes CB_*(M_n^c(k))$ is homotopic to the diagonal
cosimplicial module $diag(CB_*(C)\otimes CB_*(M_n^c(k)))$ which is
$CB_*(M_n^c(C))$.  We also know from Theorem~\ref{BarHomotopy} that
$CB_*(M_n^c(k))$ has a $gl_n^c(k) = Lie(M_n^c(k))$--comodule structure
and this coaction is null-homotopic. Corollary~\ref{Coreductivity}
implies that the differential graded $k$--module
$CB_*(M_n^c(k))^{gl_n^c(k)}$ is a direct summand of $CB_*(M_n^c(k))$ as
a differential graded $gl_n^c(k)$--comodules.  Then we can conclude that
$CB_*(M_n^c(k))^{gl_n^c(k)}$ is homotopic to the complex
$CB_*(M_n^c(k))$ since both complexes are contractible.  By using
Corollary~\ref{DualIsomorphism} we conclude
\begin{align}
CB_*(M_n^c(C))\simeq diag\left(CB_*(C)\otimes J_n^*\right)
\end{align}
where $J_n^*$ denotes the differential graded $k$--module
$\text{Hom}_k(CB_*(M_n^a(k))_{gl_n(k)},k)$.  Now, according to
\cite[Chapter 9]{Loday:CyclicHomology} Weyl's invariant theory provides
us
\begin{align}
\left(M^a_n(k)^{\otimes m}\right)_{gl_n(k)} \cong k[\Sigma_m^{ad}]
\end{align}
for $n\geq m$ where $k[\Sigma_m^{ad}]$ denotes the $\Sigma_m$--module 
$k[\Sigma_m]$ with the adjoint action.  This immediately implies that 
for $n\geq m$ one has
\begin{align}
J^m_n \cong k[\Sigma_m^{ad}]^\vee \cong k[\Sigma_m^{ad}]
\end{align}
since $k[\Sigma_m^{ad}]$ is a finite $\Sigma_m$--module.  

Consider the projective system of coalgebras
$M_{n+1}^c(k)\xra{p_{n+1}}M_n^c(k)$ where 
\begin{align}
p_n(e_{ij})
  = \begin{cases}
    e_{ij}  &  \text{ if } 1\leq i,j\leq n\\
    0       &  \text{ if } i=n+1 \text{ or } j=n+1
    \end{cases}
\end{align}
which induces a projective system of differential graded modules
$CB_*(M_{n+1}^c(k))\xra{p_{n+1}}CB_*(M_n^c(k))$ and also
$J^*_{n+1}:=CB_*(M_{n+1}^c(k))^{gl_{n+1}^c(k)}\xra{p_{n+1}}
CB_*(M_n^c(k))^{gl_n^c(k)}=:J^*_n$.  Since each $p_*$ is surjective, the
system satisfies Mittag--Leffler condition \cite[pg
82]{Weibel:HomologicalAlgebra}.  Thus
\begin{align}
0 = & \lim_n H_*(J^*_n)          \cong H_*(\lim_n J^*_n)
\end{align}
The contractible differential graded algebra $CB_*(M_n^c(C))$ is
homotopic to the diagonal cosimplicial module $diag(CB_*(C)\otimes
J_n^*)$ which in turn is homotopic to $diag(CB_*(C)\otimes J^*)$ where
$J^*=\lim_n J^*_n$.  The morphisms $J^*_{n+1}\xra{p_{n+1}}J^*_n$ are
easy to describe:
\begin{align}
p_{n+1}(\sigma)
 = \begin{cases}
     \sigma  & \text{ if } \sigma\in\Sigma_n\subseteq\Sigma_{n+1}\\
     0       & \text{ otherwise}
     \end{cases}
\end{align}
which means $J^* \cong k[\Sigma_*^{ad}]$.  Thus, an element $\Psi$
from $diag(CB_m(C)\otimes J^m)$ is of the form $\sum_i
(c^1_i\otimes\cdots\otimes c^m_i\otimes \sigma_i)$ for some $c^j_i$ from
$C$ and $\sigma_i$ from $\Sigma_m^{ad}$.  Then one can define a
coassociative cocommutative coproduct
\begin{align}
\sum_i\sum_{P\sqcup Q={\bf n},\ \sigma_i(P)=P}
   ({\bf c}^P_i\otimes\sigma^P_i)\otimes({\bf c}^Q_i\otimes\sigma^Q_i)
\end{align}
where ${\bf n}=\{1<\cdots<n\}$ and the sum is taken over all partitions
$P\sqcup Q ={\bf n}$ of ${\bf n}$ such that $\sigma_i(P)=P$.  The
element $\sigma^P_i$ denotes $\sigma_i$ considered as a set endomorphism
of ${\bf n}$ restricted to the subset $P$.  For a subset
$P=\{p_1<\cdots<p_m\}$ of ${\bf n}$ the element ${\bf c}^P_i$ is defined
to be $c^{p_1}_i\otimes\cdots c^{p_m}_i$.  We also define a product
structure 
\begin{align}
({\bf c}\otimes\sigma)\cdot({\bf d}\otimes\delta)
 = ({\bf c}\otimes{\bf d}\otimes\sigma\oplus\delta)
\end{align}
where $({\bf c}\otimes\sigma)$ and $({\bf d}\otimes\delta)$ are
arbitrary elements from $diag(CB_*(C)\otimes k[\Sigma_*^{ad}]$.  We
define $\sigma\oplus\delta$ as
\begin{align}
(\sigma\oplus\delta)(i)
 = \begin{cases}
   \sigma(i)     & \text{ if } 1\leq i\leq p\\
   \delta(i-p)+p & \text{ if } p+1\leq i\leq p+q
   \end{cases}
\end{align}
where $\sigma$ is from $\Sigma_p$ and $\delta$ is from $\Sigma_q$.

Passing to the anti-symmetrization of $diag(CB_*(C)\otimes
k[\Sigma_*^{ad}])$, one considers the collection
\begin{align}
\epsilon_m (C^{\otimes m}\otimes k[\Sigma_m^{ad}])
\cong {}_{\Sigma_m} (C^{\otimes m}\otimes k[\Sigma_m^{ad}])
\cong C^{\otimes m}\eotimes{\Sigma_m} k[\Sigma_m^{ad}]
\end{align}
for $m\geq 1$.  We must prove that the coproduct we defined above
extends to the anti-symmetrization.  For this purpose consider
$\delta\in \Sigma_m$ and ${\bf c}=(c^1\otimes\cdots\otimes c^m)$ from
$CB_m(C)$.  Then ${\bf c}\cdot\delta = c^{\delta(1)}\otimes\cdots
\otimes c^{\delta(n)}$ and
\begin{align}
\Delta({\bf c}\cdot\delta\otimes\sigma)
 = & \sum_{P\sqcup Q={\bf n},\ \sigma(P)=P}
          (({\bf c}\cdot\delta)^P\otimes\sigma^P)\otimes
          ((\delta\cdot{\bf c})^Q\otimes\sigma^Q)\\
 = & \sum_{\delta(P)\sqcup\delta(Q)=\delta({\bf n}),
           \ {}^\delta\sigma(\delta(P))=\delta(P)}
          ({\bf c}^{\delta(P)}\otimes{}(^\delta\sigma)^{\delta(P)})\otimes
          ({\bf c}^{\delta(Q)}\otimes{}(^\delta\sigma)^{\delta(Q)})\\
 = & \sum_{U\sqcup V={\bf n},\ {}^\delta\sigma(U)=U}
          ({\bf c}^U\otimes({}^\delta\sigma)^{U})\otimes
          ({\bf c}^V\otimes({}^\delta\sigma)^{V})\\
 = & \Delta({\bf c}\otimes {}^\delta\sigma)
\end{align}
which proves that the comultiplication extends to the new setting.  It
is easy to see that the multiplication we defined earlier extends to
the anti-symmetrization too.

The differential graded algebra $diag(CB_*(C)\eotimes{\Sigma_*}
k[\Sigma_*^{ad}])$ which consists of $\{C^{\otimes m}\eotimes{\Sigma_m}
k[\Sigma_m]\}_{m\geq 1}$ is homotopic to $\lim_n
CE^\text{sym}_*(gl_n^c(C))$.  Note that the induced projective system of
differential graded $k$--modules
$CE^\text{sym}_*(gl_{n+1}^c(C))\xra{p_{n+1}} CE^\text{sym}_*(gl_n^c(C))$
still satisfies Mittag--Leffler condition.  Define $gl^c(C)=\lim_n
gl_n^c(C)$ and observe that
\begin{align}
\lim_n CE^\text{sym}_*(gl_n^c(C))
  \cong CE^\text{sym}_*(\lim_n gl_n^c(C)) 
  \simeq diag(CB_*(C)\eotimes{\Sigma_*} k[\Sigma_*^{ad}])
\end{align}
According to \cite[10.2.15]{Loday:CyclicHomology} the collection
$\{C^{\otimes m}\otimes k[\Sigma_m^{ad}]\}_{m\geq 1}$ is a cocommutative
Hopf algebra with respect to the coproduct and product structures
defined above.  Then $CE^\text{sym}_*(gl^c(C)) \simeq \{C^{\otimes
m}\eotimes{\Sigma_m}k[\Sigma_m^{ad}]\}_{m\geq 1}$ is a graded
commutative/cocommutative Hopf algebra.  According to
\cite{MilnorMoore:HopfAlgebras} the Hopf algebra
$CE^\text{sym}_*(gl^c(C))$ is generated as an exterior algebra over the
submodule of primitive elements.  From the definition of the coproduct
one can see that elements of the form $\sum_i ({\bf
c}_i\otimes\sigma_i)$ are primitive where $\sigma_i$ is in the conjugacy
class of the element $\tau_m$ and ${\bf c}$ is in $C^{\otimes m}$.  With
this description, one immediately sees that the submodule of primitive
elements is also a differential graded submodule.  If $Cl(\tau_m)$
denotes the conjugacy class of the element $\tau_m$ and $C_{\tau_m}$
denotes the subgroup of centralizers of $\tau_m$ in $\Sigma_m$.  Then we
have
\begin{align}
Prim_m CE^\text{sym}_*(gl^c(C))
 =  C^{\otimes m}\eotimes{\Sigma_m}k[Cl(\tau_m)]
 \cong C^{\otimes m}\eotimes{\Sigma_m}k[\Sigma_m/C_{\tau_m}]
 \cong C^{\otimes m}\eotimes{\Sigma_m}k[\Sigma_m]\eotimes{C_{\tau_m}}k
\end{align}
which is a differential graded submodule of $CE^\text{sym}_*(gl^c(C))$.
However, the subgroup of centralizers of $\tau_m$ is $C_m$ the cyclic
group generated by $\tau_m$.  Therefore
\begin{align}
Prim_m CE^\text{sym}_*(gl^c(C))
 \cong C^{\otimes m}\eotimes{C_m}k \cong CC^\lambda_m(C)[+1]
\end{align}
as differential graded modules.  Moreover, since the submodule of
primitive elements is a differential graded submodule and since $k$ is a
field of characteristic 0 we get
\begin{align}
H^\text{Lie}_*(gl^c(C))
:= H_*CE^\text{sym}_*(gl^c(C))
\cong H_* \Lambda^*CC^\lambda_*(C)[+1]
\cong \Lambda^* HC_*(C)[+1]
\end{align}

\bibliographystyle{plain} 
\bibliography{bibliography}

\vspace{2cm}
\noindent{\small\sc 
Department of Mathematics, Ohio State University,
Columbus, Ohio 43210, USA}

\noindent{\it E-mail address:}\ {\tt kaygun@math.ohio-state.edu}

\end{document}